\documentclass[11pt, twoside, eqno]{article}
\usepackage[margin=3cm]{geometry}
\usepackage{latexsym}
\usepackage{graphicx}
\usepackage{amsmath,amsthm,amssymb,amsfonts}

\usepackage{xpatch}
\makeatletter
\xpatchcmd{\@thm}{\thm@headpunct{.}}{\thm@headpunct{}}{}{}

\renewcommand\@biblabel[1]{#1.}

\makeatother

\theoremstyle{definition}

\theoremstyle{remark}

\numberwithin{equation}{section}

\begin{document}
\noindent{\bf\Large Comment to: ``Generalized hyperideals in locally 
associative left almost semihypergroups"}\bigskip

\noindent{\bf\large Niovi Kehayopulu}\bigskip

This paper has been submitted to New York Journal of Mathematics on 
Sun, January 31, 2021.

\bigskip

\bigskip

\begin{abstract}This note is written to show that the definition of 
the ${\cal L}{\cal A}$-hypergroupoids in [5] should be corrected and 
that it is not enough to replace the multiplication ``$\cdot$" of an 
${\cal L}{\cal A}$-groupoid by the hyperoperation ``$\circ$" to pass 
from an ${\cal L}{\cal A}$-groupoid to an ${\cal L}{\cal 
A}$-hypergroupoid.
\bigskip

\end{abstract}

According to the introduction, ``a left almost semigroup (${\cal 
L}{\cal A}$-semigroup) is a groupoid $S$ whose elements satisfy the 
following left inversive law $(ab)c=(cb)a$ for all $a,b,c\in S$".
This structure being a groupoid and not a semigroup, the term
${\cal L}{\cal A}$-semigroup is wrong and it should be replaced by 
${\cal L}{\cal A}$-groupoid also called Abel-Grassmann's groupoid.
Some authors use the concept ``Non associative ordered semigroup" 
that is certainly wrong as an ordered semigroup cannot be non 
associative. Again by the introduction, an ${\cal L}{\cal 
A}$-groupoid is a structure between a groupoid and a commutative 
semigroup. In what sense\,? And an ${\cal L}{\cal A}$-groupoid is a 
generalization of semigroup. How is it possible\,?

For a nonempty set $H$, denote by ${\cal P}^*(H)$ the set of nonempty 
subsets of $H$.

An hypergroupoid is a nonempty set $H$ with a mapping $$\circ : 
H\times H \rightarrow {\cal P}^*(H) \mid (a,b)\rightarrow a\circ b$$ 
(called hyperoperation as it assigns to each couple $(a,b)$ of 
elements of $H$ a subset -instead of an element- of $H$).

For two nonempty subsets $A,B$ of $H$, the authors denote
$$A\circ B=\bigcup\limits_{a \in A,\,b \in B} {a\circ b}, \;\; a\circ 
A=\{a\}\circ A \mbox { and } a\circ B=\{a\}\circ B. $$Then, they 
define the ${\cal L}{\cal A}$-semihypergroup as an hypergroupoid such 
that $$(x\circ y)\circ z=(z\circ y)\circ z \mbox { for all } x,y,z\in 
H.$$
First of all, we cannot use the term ``${\cal L}{\cal 
A}$-semihypergroup" (left almost semihypergroup) exactly as we cannot 
use the term ${\cal L}{\cal A}$-semigroup. ${\cal L}{\cal 
A}$-hypergroupoid is the correct one and this is what we will use in 
the present note.

We cannot write $A\circ B$, as ``$\circ$" is an ``operation" between 
elements. We cannot use the same symbol both for elements and sets; 
if we do that, a great confusion erases (see also [2]).

We cannot write $(x\circ y)\circ z$ as $x\circ y$ is a subset of $H$, 
$z$ an element of $H$ and ``$\circ$" is an ``operation" between 
elements. Even if we identify the $z$ by the singleton $\{z\}$, we 
cannot write $(x\circ y)\circ \{z\}$ for the reason mentioned above.

It is very difficult (not to say impossible) to check the examples of 
${\cal L}{\cal A}$-hypergroupoids by hand, so we should know the way 
they have been constructed. Or should be given a method to check 
their validity. There is nothing like that in the bibliography.

According to Example 4,
the set $H=\{a,b,c,d,e\}$ with the hyperoperation $\circ$ given by 
Table 1, is an ${\cal L}{\cal A}$-hypergroupoid.
\begin{center}
$\begin{array}{*{20}{c}}
{\circ}&\vline& a&\vline& b&\vline& c&\vline& d&\vline& e\\
\hline
a&\vline& {\{ a\} }&\vline& {\{ a\} }&\vline& {\{ a\} }&\vline& {\{ 
a\} }&\vline& {\{ a\} }\\
\hline
b&\vline& {\{ a\} }&\vline& {\{ a,e\} }&\vline& {\{ a,e\} }&\vline& 
{\{ a,c\} }&\vline& {\{ a,e\} }\\
\hline
c&\vline& {\{ a\} }&\vline& {\{ a,e\} }&\vline& {\{ a,e\} }&\vline& 
{\{ a,b\} }&\vline& {\{ a,e\} }\\
\hline
d&\vline& {\{ a\} }&\vline& {\{ b\} }&\vline& {\{ c\} }&\vline& {\{ 
d\} }&\vline& {\{ e\} }\\
\hline
e&\vline& {\{ a\} }&\vline& {\{ a,e\} }&\vline& {\{ a,e\} }&\vline& 
{\{ a,e\} }&\vline& {\{ a,e\} }
\end{array}$\bigskip

Table 1
\end{center}We wrote $\{a\}$ in the table instead of the $a$ in [5].
According to [5], to show that this is an ${\cal L}{\cal 
A}$-hypergroupoid we have to show that $(x\circ y)\circ z=(z\circ 
y)\circ x$ for all $x,y,z\in H$; as so the $(a\circ b)\circ d$, for 
example, should be equal to $(d\circ b)\circ a$. We have $(a\circ 
b)\circ d=\{a\}\circ d$ and $(d\circ b)\circ a=\{b\}\circ a$, while 
both the
$\{a\}\circ d$ and the $\{b\}\circ a$ are without meaning.

Throughout the paper strange symbols like
$$[(H\circ R)\circ (H\circ R^{(m+n-2)})]\circ (H\circ L)\circ (H\circ 
L^ {(m+n-2)})]$$  have been used (see, for example, p. 1068, l. 5). 
Symbols without any sense.

What is the $A^n$\,? We are not in a semigroup where $A^{n}=A\,\cdot 
A\cdot \cdot \cdot\, A$. In an ${\cal L}{\cal A}$-hypergroupoid can 
we write $A*A* \cdot \cdot \cdot *A$ without using parentheses? 
Certainly not. So, for an ${\cal L}{\cal A}$-hypergroupoid the 
concept of the $(m,n)$-hyperideal, for arbitrary $m,n$ cannot be 
defined. The concept of $(m,n)$-regularity of this structure cannot 
be defined for the same reason. What is the $a^{m}$ for an arbitrary 
natural number $m$\,? \bigskip

Let us first correct the definitions in [5].

Let $H$ be a nonempty set and ``$\circ$" be an hyperoperation on $H$. 
For two nonempty subsets $A$, $B$ of $H$ denote by ``$*$" the 
operation on ${\cal P}^* (H)$ (induced by $\circ$) defined by
$$*: {\cal P}^* (H)\times {\cal P}^* (H)\rightarrow {\cal P}^* 
(H)\mid (A,B) \rightarrow A*B:= \bigcup\limits_{a \in A,\,b \in B} 
{a\circ b}.$$An ${\cal L}{\cal A}$-hypergroupoid is a nonempty set 
$H$ with an hyperoperation ``$\circ$" such that $$(a\circ 
b)*\{c\}=(c\circ b)*\{a\} \mbox { for all } a,b,c\in H.$$As one can 
easily see, for any $a,b$ we have $\{a\}*\{b\}=a\circ b$ (see, for 
example [1--3]).

With this definition, let us check if the Example 4 in [5] (given by 
Table 1 above) is correct: It is wrong as, for example,

$(d\circ d)*\{b\}=\{d\}*\{b\}=d\circ b=\{b\}$,

$(b\circ d)*\{d\}=\{a,c\}*\{d\}=\bigcup\limits_{x \in \{ a,\,c\} } 
{x\circ d}=(a\circ d)\cup (c\circ d)=\{a\}\cup \{a,b\}=\{a,b\}$ and

$(d\circ d)*\{b\}\not=(b\circ d)*\{d\}$.\medskip

According to Example 5 of the paper, the set $H=\{a,b,c,d\}$ defined 
by Table 2 is an ${\cal L}{\cal A}$-hypergroupoid.
\begin{center}
$\begin{array}{*{20}{c}}
{\circ}&\vline& a&\vline& b&\vline& c&\vline& d\\
\hline
a&\vline& {\{ a\} }&\vline& {\{ a\} }&\vline& {\{ a\} }&\vline& {\{ 
a\} }\\
\hline
b&\vline& {\{ a\} }&\vline& {\{ a,b,c,d\} }&\vline& {\{ a,b,c\} 
}&\vline& {\{ a,b,c\} }\\
\hline
c&\vline& {\{ a\} }&\vline& {\{ a,b,c\} }&\vline& {\{ a,b,c\} 
}&\vline& {\{ a,b,c\} }\\
\hline
d&\vline& {\{ a\} }&\vline& {\{ a,b,c\} }&\vline& {\{ a,b,c\} 
}&\vline& {\{ a,b,c\} }
\end{array}$\bigskip

Table 2
\end{center}
This is also wrong as

$(c\circ d)*\{b\}=\{a,b,c\}*\{b\}=\bigcup\limits_{x \in \{ a,b,c\} } 
{x\circ b}=(a\circ b)\cup (b\circ b)\cup (c\circ b)=\{a,b,c,d\}$,

$(b\circ d)*\{c\}=\{a,b,c\}*\{c\}=\bigcup\limits_{x \in \{ a,b,c\} } 
{x\circ c}=(a\circ c)\cup (b\circ c)\cup (c\circ c)=\{a,b,c\}$ but

$(c\circ d)*\{b\}\not=(b\circ d)*\{c\}$.\medskip

In what follows, the aim is to show that is not enough to pass from 
an ${\cal L}{\cal A}$-groupoid to an ${\cal L}{\cal A}$-hypergroupoid 
by replacing the multiplication ``$\cdot$" of the ${\cal L}{\cal 
A}$-groupoid by the hyperoperation ``$\circ$" of the ${\cal L}{\cal 
A}$-hypergroupoid.

The paper in [5] is the paper in [4] with the only difference that 
the multiplication ``$\cdot$" in [4] has been replaced by ``$\circ$" 
in [5].

In fact,

Lemma 1 in [5] is the Lemma 7 in [4];

Theorem 1 in [5] is the Theorem 6 in [4];

Theorem 2 in [5] is the Theorem 7 in [4];

Theorem 3 in [5] is the Theorem 8 in [4];

Lemma in [5] is the Lemma 8 in [4];

Theorem 2 in [5] is the Theorem 7 in [4];

Theorem 4 in [5] is the Theorem 9 in [4];

Corollary 1 in [5] is the Corollary 5 in [4];

Theorem 5 in [5] is the Theorem 10 in [4];

Lemma 3 in [5] is the Lemma 9 in [4];

Corollary 2 in [5] is the Corollary 6 in [4];

Theorem 6 in [5] is the Theorem 11 in [4];

Theorem 7 in [5] is the Theorem 12 in [4];

Theorem 8 in [5] is the Theorem 13 in [4].\medskip

Finally, 29 papers have been cited in References while the paper is 
based only on [4] not cited in References.

It might be mentioned here, though it is out of the scope of the 
present note, that the paper in [6] is an exact copy (word by word) 
by the paper in [4] with the only difference that 17 papers have been 
cited in References in [6] instead of 13 in [4].\medskip

As far as the fourth section (Conclusions) is concerned, it is 
obvious that an ${\cal L}{\cal A}$-groupoid is not ${\cal L}{\cal 
A}$-hypergroupoid in general (no Reference to this was needed). From 
this the authors conclude that the results of their paper generalize 
the results of ${\cal L}{\cal A}$-groupoids.
But to put ``$\circ$" instead of the multiplication ``$\cdot$" of the 
${\cal L}{\cal A}$-groupoid and use strange symbols without any 
sense, we can never say that the results of the ${\cal L}{\cal 
A}$-groupoids have been generalized. On the other hand, the results 
in [4] should be checked as well since some lemmas were needed on 
which they were based.
How we prove the results in [4] using the $A^{n}$\,? What is the 
$A^{n}$\,? According to this section, the results of the paper extend 
the result by M. Akram, N. Yaqoob and M. Khan [On $(m,n)$-ideals in 
LA-semigroups, Appl. Math. Sci. (Ruse) 7, no. 41--44 (2013), 
2187--2191] ([1] in References of [5]) and by Q. Mushtaq and S.M. 
Yusuf [On locally associative LA-semigroups, J. Natur, Sci. Math. 19, 
no. 1 (1979), 57--62] ([21] in References of [5]).
The results in [1] have nothing to do with this paper while [21] is 
an old paper for which we have only the review MR0596763 (82a:20081) 
by M. Friedberg in MathSciNet; does not seem to be related to [5].
\bigskip

\noindent (Niovi Kehayopulu) DEPARTMENT OF MATHEMATICS, UNIVERSITY OF 
ATHENS, 15784 PANEPISTIMIOPOLIS, GREECE\\
nkehayop@math.uoa.gr

\end{document}